\documentclass[12pt]{article}

\usepackage{amsmath,amsthm}

\parskip=\medskipamount
\def\a{\alpha}

\def\cinfty#1{C^{\scriptscriptstyle\infty}(#1)}
\def\vectorfields#1{{\cal X}(#1)}
\def\ov#1{\overline{#1}}
\def\G{\Gamma}

\def\del#1{{\nabla}_{#1}}

\def\fpd#1#2{\frac{\partial #1}{\partial #2}}
\def\R{{\rm I\kern-.20em R}}
\def\sode{{\sc Sode}}

\def\proof{{\sc Proof}}
\newtheorem{prop}{\bf Proposition}

\newtheorem{thm}{\bf Theorem}
\newtheorem{dfn}{\bf Definition}

\def\r{\rho}

\def\z{\zeta}

\def\base#1{\pmb{e}_{#1}}
\def\cosyst{(e_0; \{ \base \a \})}
\def\baseX#1{{\mathcal X}_{#1}}
\def\baseV#1{{\mathcal V}_{#1}}
\def\baseH#1{{\mathcal H}_{#1}}

\def\H#1{{#1}^{\scriptscriptstyle H}}
\def\V#1{{#1}^{\scriptscriptstyle V}}
\def\subH#1{{#1}_{\scriptscriptstyle H}}
\def\subV#1{{#1}_{\scriptscriptstyle V}}

\def\tpi{\tilde\pi}

\def\prol#1#2{{T^{#1}{#2}}}
\def\prolr#1{\prol{\rho}{#1}}

\newcommand{\mybox}[1]{\makebox(0,0){\footnotesize{#1}}}
\newlength{\savelen}

\begin{document}

\title{Note on generalised connections and affine bundles}

\author{T.\ Mestdag, W.\ Sarlet\\
{\small Department of Mathematical Physics and Astronomy }\\
{\small Ghent University, Krijgslaan 281, B-9000 Ghent, Belgium}\\[2mm]
E.\ Mart\'{\i}nez\\
{\small Departamento de Matem\'atica Aplicada}\\
{\small Facultad de Ciencias, Universidad de Zaragoza}\\
{\small Pedro Cerbuna 12, E-50009 Zaragoza, Spain}}

\maketitle

\begin{quote}
{\bf Abstract.} {\small We develop an alternative view on the concept of
connections over a vector bundle map, which consists of a horizontal lift
procedure to a prolonged bundle. We further focus on prolongations to an
affine bundle and introduce the concept of affineness of a generalised
connection.}
\end{quote}

%AMS-classification: ?????

\section{Introduction}

There has been a lot of interest, recently, in potential applications of Lie
algebroids in physics, control theory and other fields of applied mathematics.
Among papers which study, in particular, aspects of Lagrangian systems on Lie
algebroids, we mention
\cite{Wein,Liber,Mart,Carin,CarMar,Clem,affalg1,affalg2}. There is of course
an enormous literature on more purely mathematical aspects of Lie algebroids,
of which we cite only the standard work \cite{Mac}, and \cite{Fernandes} for
its particular relevance to this paper.

Our recent joint work in the field finds its roots in searching for the right
geometrical model for a kind of time-dependent generalisation of `Lagrangian
mechanics' on Lie algebroids. Since ordinary time-dependent mechanics is
usually described on the first-jet space $J^1M$ of a manifold $M$ fibred over
$\R$ (see for example \cite{Crampin,Mangia1}), the direct model for the
generalisation we had in mind was a kind of Lie algebroid structure, whose
anchor map takes values in $J^1M$ rather than $TM$. This was explored in
detail in \cite{affalg1}, which in turn rose interest in the more general
features of having a Lie algebroid structure on an affine bundle, without the
requirement that the base manifold be fibred over $\R$. Those ideas were
developed in \cite{affalg2} and to some extent (that is without reference to
dynamical systems) also in \cite{GraGra}.

A continuation of this work is in preparation, in particular with the purpose of
bringing a suitable theory of
connections into the picture of dynamical systems on affine algebroids. But the path to
these further developments has led us to discover some general features on connections and
affine spaces, which do not require a Lie algebroid structure and seem well worth being
brought under the attention separately. This brings us to the subject matter of the
present paper.

In Section~2, the main objective is to discuss two interesting constructions
from the recent literature on generalised connections and algebroids which,
when brought together in a unifying picture, will open the way to explain in
detail how they are related. Both constructions may have their roots in the
theory of Lie algebroids, but have been formulated recently in the more
general framework where a vector bundle has a kind of anchor map, but need not
be equipped with a Lie algebra structure on the real vector space of its
sections. The first topic we are referring to is the notion of generalised
connection on a vector bundle map, as introduced by Cantrijn and Langerock
\cite{FransBavo}, inspired by a similar construction on Lie algebroids by
Fernandes \cite{Fernandes}. The second is the idea of prolongation, which has
been discussed in the context of Lie algebroids, for example in
\cite{Mart,affalg1,affalg2}, but, as shown in \cite{affalg2}, can also be
defined without the need of a Lie algebra structure. Also relevant is work by
Popescu, who in fact already developed the same ideas in the case that all
bundles involved are vector bundles; for this we refer to \cite{Pop} and
references therein. We will arrive, in Section~2, at an alternative view on
the generalised connections of \cite{FransBavo}. But let us mention here
already that this alternative view can be developed without needing the
generalised connection idea of \cite{FransBavo}. This is in fact one of the
main discoveries of \cite{Mart} and \cite{Pop} and it is being explored to
full extent in \cite{MartII}. The purpose of the present note, however, is to
explain in detail the interrelationship between the two ideas.

In Section~3, we focus on the case of the prolongation of an
affine bundle $E\rightarrow M$ over a vector bundle $V\rightarrow
M$. We show that bringing the bidual of $E$ into the picture
enables us to give a clear and concise definition of the concept
of an affine connection over a vector bundle map and prove a
result about the equivalent characterisation of such a connection
via a kind of covariant derivative operator. The relevance of
these results for the future developments we have in mind is
briefly indicated in the final section.

\section{Connections over a vector bundle map and the horizontal subbundle of
a prolonged bundle}

We start by recalling the prolongation idea, as developed in \cite{affalg2}.

Let $\mu: P\rightarrow M$ be an arbitrary fibre bundle and $\tau:V\rightarrow M$ a vector
bundle. Assume there exists an {\em anchor map} $\rho: V \rightarrow TM$, which for the time
being is just a vector bundle morphism.

\begin{dfn} The $\rho$-prolongation of $\mu:P\rightarrow M$ is the bundle $\mu^1: \prolr{P} \rightarrow
P$, constructed as follows: (i) the total space $\prolr{P}$ is the
total space of the pullback bundle $\rho^*TP$
\begin{equation}
\prolr{P} = \{(v,X_p) \in V \times TP \mid \, \rho(v) = T\mu(X_p) \} \label{LVP};
\end{equation}
(ii) if $\rho^1$ denotes the projection of $\rho^*TP$ into $TP$
and $\tau_P$ is the tangent bundle projection, then
$\mu^1=\tau_P\circ\rho^1$.
\end{dfn}

The situation is summarised in the following diagram, whereby the
projection on the first element of a pair $(v,X_p)\in \prolr{P}$
is denoted by $\mu^2$.

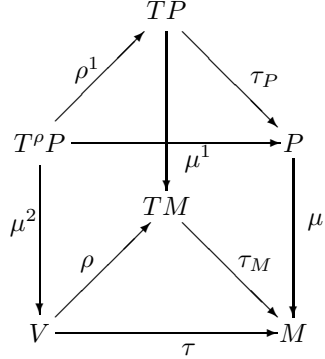
\setlength{\unitlength}{12pt}

\begin{picture}(40,13)(-6,2)

\put(7.5,9.5){\vector(1,0){6.6}}
\put(11,13){\vector(1,-1){3}}
\put(7,10.2){\vector(1,1){2.8}}

\put(6.5,9.5){\mybox{$\prolr{P}$}} \put(14.5,9.5){\mybox{$P$}}
\put(10.5,13.7){\mybox{$TP$}}

\put(7,3.5){\vector(1,0){7}}
\put(11,7){\vector(1,-1){3}}
\put(7,4){\vector(1,1){3}}

\put(6.5,3.5){\mybox{$V$}} \put(14.5,3.5){\mybox{$M$}} \put(10.5,7.5){\mybox{$TM$}}

\put(6.5,8.8){\vector(0,-1){4.7}}
\put(10.5,13){\vector(0,-1){5}}
\put(14.5,9){\vector(0,-1){5}}

\put(11.2,3){\mybox{$\tau$}} \put(8.0,5.7){\mybox{$\rho$}}
\put(13.3,5.7){\mybox{$\tau_M$}}

\put(11.5,9){\mybox{$\mu^1$}} \put(8.0,11.8){\mybox{$\rho^1$}}
\put(13.6,11.5){\mybox{$\tau_P$}}

\put(6.0,7){\mybox{$\mu^2$}} \put(15.2,7){\mybox{$\mu$}}

\end{picture}
\setlength{\unitlength}{\savelen}

One can think of the bundle $\mu^1:\prolr{P}\rightarrow P$ as a kind of
generalisation of a tangent bundle. Obviously, the standard tangent bundle
fits into the picture: it suffices to take $V=TM$ and $\rho=id_{TM}$. More
interestingly, if we have two bundles $\mu_i:P_i\rightarrow M$ ($i=1,2$), and
a bundle map $f$ (over the identity on $M$) between them, then the tangent map
$Tf:TP_1\rightarrow TP_2$ extends to a map $\prolr{P_1} \rightarrow
\prolr{P_2}: (v, X_p) \mapsto (v, Tf(X_p))$. Indeed, we have
$T\mu_2(Tf(X_p))=T\mu_1(X_p)=\rho(v)$. There is more to say about the
tangent-bundle-like behaviour of $\prolr{P}$, but we will not elaborate on
that here.

Coming back to the diagram above, an element of $\prolr{P}$ is called {\em
vertical\/} if it is in the kernel of the projection $\mu^2$. The set of all
vertical elements in $\prolr{P}$ is a vector subbundle of $\mu^1$ and will be
denoted by ${\mathcal V}^\rho P$. If $(0,Q)\in {\mathcal V}^\rho P$, then
$Q=\rho^1(0,Q)$ will also be vertical in $TP$, since $T\mu(Q)=\rho(0)=0$. The
idea of arriving at a notion of horizontality on $TP$, adapted to the presence
of the anchor map in the picture, lies at the basis of the following concept,
introduced in \cite{FransBavo}.

\begin{dfn} A $\rho$-connection on $\mu$ is a
linear bundle map $h: \mu^* V \rightarrow TP$ (over the identity on $P$),
such that $\rho \circ p_V = T\mu \circ h$, where $p_V$ is the projection of $\mu^*V$ onto $V$.
\end{dfn}

There is a quite striking similarity between our first diagram and the one we
can draw here for the illustration of all spaces involved in the definition of
a $\rho$-connection:

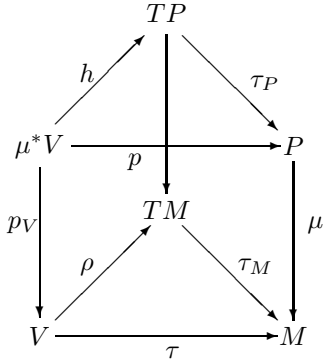
\setlength{\unitlength}{12pt}

\begin{picture}(40,13)(-6,2)

\put(7.5,9.5){\vector(1,0){6.6}} \put(11,13){\vector(1,-1){3}}
\put(7,10.2){\vector(1,1){2.8}}

\put(6.5,9.5){\mybox{$\mu^*V$}} \put(14.5,9.5){\mybox{$P$}}
\put(10.5,13.7){\mybox{$TP$}}

\put(7,3.5){\vector(1,0){7}} \put(11,7){\vector(1,-1){3}}
\put(7,4){\vector(1,1){3}}

\put(6.5,3.5){\mybox{$V$}} \put(14.5,3.5){\mybox{$M$}}
\put(10.5,7.5){\mybox{$TM$}}

\put(6.5,8.8){\vector(0,-1){4.7}} \put(10.5,13){\vector(0,-1){5}}
\put(14.5,9){\vector(0,-1){5}}

\put(10.7,3){\mybox{$\tau$}} \put(8.0,5.7){\mybox{$\rho$}}
\put(13.3,5.7){\mybox{$\tau_M$}}

\put(9.5,9){\mybox{$p$}} \put(8.0,11.8){\mybox{$h$}}
\put(13.6,11.5){\mybox{$\tau_P$}}

\put(6.0,7){\mybox{$p_V$}} \put(15.2,7){\mybox{$\mu$}}

\end{picture}
\setlength{\unitlength}{\savelen}

Note that points in the image $\rho^1(\prolr{P})$ can be vertical in $TP$ when
the corresponding point in the domain is not vertical in $\prolr{P}$ (because
$\rho$ need not be injective). This is related to the observation that ${\rm
Im}\,h$ can have a non-empty intersection with the vertical vectors on $P$. As
discussed in detail in \cite{FransBavo}, ${\rm Im}\,h$ will in general also
fail to determine a full complement to the vertical vectors on $P$. That is
why one refers to a $\rho$-connection on $\mu$ also as a `generalised
connection'.

The point we would like to emphasise, however, is that it is
perhaps not such a good idea to concentrate on horizontality on
$TP$. Instead, as one may conjecture from an inspection of the two
diagrams above, the better fibration to look for horizontality in
this framework is the prolonged bundle $\mu^1:\prolr{P}\rightarrow
P$. In other words, we think it is important to bring the pictures
of {\em $\rho$-prolongation} and {\em $\rho$-connection} together
into the following scheme.

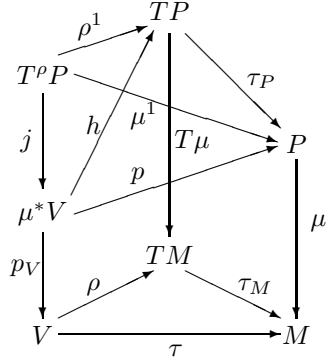
\setlength{\unitlength}{12pt}

\begin{picture}(40,13)(-6,2)

\put(7.5,7.3){\vector(3,1){6.4}} \put(7.5,11.7){\vector(3,-1){6.4}}
\put(11,13){\vector(1,-1){3}} \put(7,12.3){\vector(3,1){2.8}}
\put(7.4,7.9){\vector(1,2){2.6}}

\put(6.5,7.3){\mybox{$\mu^*V$}} \put(6.5,11.7){\mybox{$\prolr{P}$}}
\put(14.5,9.5){\mybox{$P$}} \put(10.5,13.7){\mybox{$TP$}}

\put(7,3.5){\vector(1,0){7}} \put(11,5.5){\vector(2,-1){3}}
\put(7,4){\vector(2,1){3}}

\put(6.5,3.5){\mybox{$V$}} \put(14.5,3.5){\mybox{$M$}}
\put(10.5,6){\mybox{$TM$}}

\put(6.5,6.7){\vector(0,-1){2.7}} \put(10.5,13){\vector(0,-1){6.4}}
\put(14.5,9){\vector(0,-1){5}} \put(6.5,11.1){\vector(0,-1){3.0}}

\put(10.7,3){\mybox{$\tau$}} \put(8.1,5.0){\mybox{$\rho$}}
\put(13.2,5.0){\mybox{$\tau_M$}}

\put(9.5,8.5){\mybox{$p$}} \put(8.0,13.2){\mybox{$\rho^1$}}
\put(13.4,11.5){\mybox{$\tau_P$}} \put(9.7,10.4){\mybox{$\mu^1$}}

\put(6.0,5.6){\mybox{$p_V$}} \put(15.2,7){\mybox{$\mu$}}
\put(6.0,9.6){\mybox{$j$}} \put(11.2,9.6){\mybox{$T\mu$}}
\put(8.1,10.1){\mybox{$h$}}

\end{picture}
\setlength{\unitlength}{\savelen}

What we propose to discuss in detail now is that, given a
$\rho$-connection on $\mu$, there is an associated, genuine
decomposition of the bundle $\mu^1$, i.e.\ a `horizontal
subspace', at each point $p\in P$, of the fibre of $\prolr{P}$,
which is complementary to the vertical subspace at $p$. In other
words, instead of considering a horizontal lift operation from
sections of $\tau$ to sections of $\tau_P$, as is done in
\cite{FransBavo}, it is more appropriate to focus on a horizontal
lift from sections of $\tau$, and by extension sections of the
pullback bundle $p$, to sections of the bundle $\mu^1$.

The fibre linear map $j: \prolr{P} \rightarrow \mu^* V: (v,Q) \mapsto
(\tau_P(Q),v)$ is surjective and its kernel is ${\mathcal V}^\rho P$.
Therefore we have the following short exact sequence of vector bundles:
\begin{equation}
0 \rightarrow {\mathcal V}^\rho P \stackrel{}{\rightarrow}
\prolr{P} \stackrel{j}{\rightarrow} \mu^* V \rightarrow 0,
\label{shortexactgeneral}
\end{equation}
where the second arrow is the natural injection.

\begin{thm} The existence of a $\rho$-connection on $\mu$ is equivalent to
the existence of a splitting $\H{}$ of the short exact sequence
(\ref{shortexactgeneral}); we have $\rho^1\circ \H{}=h$.
\end{thm}

\proof \ Let $h: \mu^* V \rightarrow TP$ be given and satisfy the requirements
of a $\rho$-connection on $\mu$. To define the `horizontal lift' of a point
$(p,v)\in\mu^*V$, as a point in $\prolr{P}$, it suffices to fix the
projections of $\H{(p,v)}$ under $\rho^1$ and $\mu^2$ in a consistent way. We
put:
\begin{equation}
\rho^1\big((p,v\H)\big) := h(p,v) \qquad \mbox{and} \qquad
\mu^2\big( (p,v \H) \big) := v. \label{Hlift}
\end{equation}
This determines effectively an element of $\prolr{P}$ since
$\rho\circ\mu^2 ((p,v \H)) = \rho(v) = \rho \circ p_V ((p,v)) =
T\mu\circ h ((p, v)) = T\mu\circ\rho^1((p,v \H))$. The horizontal
lift is obviously a splitting of (\ref{shortexactgeneral}), since
by construction $j\left((p,v
\H)\right)=\left(\tau_P(h(p,v)),v\right)=(p,v)$.

Conversely, if a splitting $\H{}$ of (\ref{shortexactgeneral}) is given, we
define $h:\mu^*V\rightarrow TP$ by $h(p,v)=\rho^1(\H{(p,v)})$. It satisfies
the required properties, i.e.\ $h$ is a linear bundle map and we have
\[
T\mu \circ h = T\mu\circ \rho^1 \circ \H{} = \rho \circ \mu^2
\circ \H{} = \rho \circ p_V\circ j\circ \H{}=\rho\circ p_V,
\]
which concludes the proof. \qed

Denoting the subbundle of $\prolr{P}$ which is complementary to
${\mathcal V}^\rho P$ by ${\mathcal H}^\rho P$, it follows that
\begin{equation}
\prolr{P}={\mathcal H}^\rho P \oplus {\mathcal V}^\rho P.
\label{directsum}
\end{equation}
An equivalent way of expressing this decomposition (analogous to what is
familiar for the case of a classical Ehresmann connection) is the following:
there exist two complementary projection operators $\subH P$ and $\subV P$ on
$\prolr{P}$, i.e.\ we have $\subH{P} + \subV{P}=id$, and
\[
{\subH{P}}^2=\subH{P}, \qquad {\subV{P}}^2=\subV{P}, \qquad
\subH{P}\circ\subV{P}=\subV{P}\circ\subH{P}=0.
\]

As usual, (\ref{shortexactgeneral}) leads to an associated short exact
sequence for the set of sections of these spaces, regarded as bundles over
$P$:
\begin{equation}
0 \rightarrow Ver(\mu^1) \stackrel{}{\rightarrow} Sec(\mu^1)
\stackrel{j}{\rightarrow} Sec(p) \rightarrow 0, \label{shortexactsec}
\end{equation}
where $Ver(\mu^1)$ denotes the set of vertical sections of $\mu^1$. The same
symbol $j$ is used for this second interpretation, so that for ${\mathcal
Z}\in Sec(\mu^1)$ and $p\in P$: $j({\mathcal Z})(p)=j({\mathcal Z}(p))$. Via
the composition with $\mu$, sections of $\tau$ can be regarded as maps from
$P$ to $V$ and, as such, are ({\em basic}) sections of $p: \mu^* V \rightarrow
P$. We will use the notations $\subV{P}$ and $\subH{P}$ also when we regard
these projectors as acting on sections of $\mu^1$, rather than points in
$\prolr{P}$.

Apart from the already mentioned applications to Lie algebroids
\cite{Fernandes,Pop}, it has recently been shown that $\rho$-connections can
be an important tool in, for example, nonholonomic mechanics \cite{BavoI},
sub-Riemannian geometry \cite{BavoII}, Poisson geometry \cite{Fernandes2} and
in control theory \cite{BavoIII}.

\subsection*{The case of linear $\rho$-connections}

Assume now that $\mu: P \rightarrow M$ now is a {\em vector\/} bundle.
Linearity of a connection is characterised in \cite{FransBavo} by an
invariance property of the map $h$ under the flow of the dilation field on
$P$. A more direct characterisation of linearity is the following. Let
$\Sigma_\lambda: P \times_M P \rightarrow P$ denote the linear combination
map: $\Sigma_\lambda(p_1,p_2)=p_1+\lambda p_2$. A $\rho$-connection on $\mu$
is said to be {\em linear} if the map $h:\mu^*V \rightarrow TP$ has the
property
\begin{equation}
h(p_1+\lambda p_2,v)= T_{(p_1,p_2)}\Sigma_\lambda\,\big(h(p_1,v),h(p_2,v)\big)
\label{linearcon},
\end{equation}
for all $(p_1,p_2)\in P\times_M P$, $\lambda\in\R$ and $v\in V$.

As is shown in \cite{FransBavo}, any operator $\nabla: Sec(\tau) \times
Sec(\mu) \rightarrow Sec(\mu)$ which is $\R$-bilinear and has the properties
\begin{equation}
\nabla_{f\zeta}\sigma = f \nabla_\zeta \sigma, \qquad  \qquad \nabla_\zeta
(f\sigma) = f\nabla_\zeta \sigma + \rho(\zeta)(f)\sigma, \label{covderiv}
\end{equation}
for all $\zeta\in Sec(\tau)$, $\sigma\in Sec(\mu)$ and $f\in\cinfty{M}$,
defines a unique {\em linear\/} $\rho$-connection on $\mu$. As usual, the
linearity of the covariant derivative operator $\nabla$ in its first argument,
implies that the value of $\nabla_\zeta \sigma$ at a point $m\in M$, only
depends on the value $\zeta$ at $m$ and thus gives rise to an operator
$\nabla_{v}: Sec(\mu) \rightarrow P_{\tau(v)}$, for each $v\in V$, determined
by
\[
\nabla_{v}{\eta}:=\nabla_{\zeta}{\eta}(m), \qquad \mbox{with} \quad \z(m)=v.
\]
In order to come to a covariant derivative along curves and a rule of parallel
transport, we make the following preliminary observation. Going back to the
overall diagram, we see two ways to go from $\prolr{P}$ to $TP$, namely the
direct map $\rho^1$ and $h\circ j$. By definition, the image for both maps
projects under $T\mu$ onto the same $\rho(v)$, so that the difference is a
vertical vector at some point $p\in P$ which, when $P$ is a vector bundle,
can be identified with an element of $P_{\mu(p)}$. With these identifications
understood, we eventually get a map from $\prolr{P}$ to $P$
which is called the connection map in \cite{FransBavo} (by analogy with the
connection map in \cite{Vilms}). Let us summarise this
by writing simply
\begin{equation}
K:=\rho^1 - h\circ j: \prolr{P}\rightarrow P \label{connnectionmap}
\end{equation}
(read: $K$ is $\rho^1 - h\circ j$, when regarded as map from $\prolr{P}$ into
$P$). The following side observation is worth being made here. In the
alternative concept of $\rho$-connections, as established by Theorem~1, it is
clear that the connection map $K$ is nothing but the vertical projector
$\subV{P}$, with a similar identification being understood (to be precise: the
isomorphism between ${\mathcal V}^\rho_p P$ and $V_pP$, followed by the
identification with $P_{\mu(p)}$ again). In fact this illustrates that the
alternative view is superior to the one expressed by Definition~2, in the
following sense. Once the importance of the space $\prolr{P}$ is recognised,
one can (in the present case that $P$ is a vector bundle) define a vertical
lift operation from $P_{\mu(p)}$ to ${\mathcal V}^\rho_p P$ in the usual way
(see the next section for more details); it extends to sections of bundles
over $P$, i.e.\ yields a vertical lift from sections of $\mu^*P\rightarrow P$
to $Sec(\mu^1)$. So, it is a matter of developing first these tangent bundle
like features of the $\rho$-prolongation, after which all tools are available
to discuss $\rho$-connections without ever needing the map $h$. This is the
main merit of the approach taken in \cite{Pop} and \cite{MartII}. For the sake
of further unifying both pictures, however, we will continue here to take
advantage of the insight which is being offered by our overall diagram.

Let now $c:I\rightarrow V$ be a $\rho$-admissible curve, which means that
$\dot{c}_M=\rho\circ c$, where $c_M=\tau\circ c$ is the projected curve in
$M$. Consider further a curve $\psi:I\rightarrow P$ in $P$ which projects on
$c_M$, i.e.\ such that $\psi_M:= \mu\circ\psi=c_M$. It follows that
$T\mu\circ\dot{\psi}=\rho\circ c$, so that such a $\psi$ actually gives rise
to a curve in $\prolr{P}$: $t\mapsto (c(t),\dot{\psi}(t))$. As a result,
making use of the map $K$, we can obtain a new curve in $P$, which is denoted
by $\nabla_c\psi$:
\begin{equation}
\nabla_c\psi(t):= K((c(t),\dot{\psi}(t)))=\dot{\psi}(t) - h((\psi(t),c(t))),
\label{covderivcurve}
\end{equation}
(the identification of $P$ with $VP$ being understood). If $\eta$
is a section of $\mu$ and $c$ is an admissible curve, then
denoting by $\psi$ the restriction of $\eta$ to that curve,
$\psi(t)=\eta(c_M(t))$, one can show that
\begin{equation}
\nabla_c\psi(t)= \nabla_{c(t)}\eta. \label{covderivcurve2}
\end{equation}
As can be readily seen from (\ref{covderivcurve}), given an admissible curve
$c$ and a point $p\in P$, finding a curve $\psi$ in $P$ which starts at $p$
and makes $\nabla_c\psi=0$ is a well-posed initial value problem for a
first-order ordinary differential equation, and hence gives rise to a unique
solution. The solution is called the {\em horizontal lift\/} of $c$ through
$p$, denoted by $c^h$. Hence, we have
\begin{equation}
\nabla_c c^h=0, \label{h-lift}
\end{equation}
and points in the image of $c^h$ are said to be obtained from $p$ by {\em
parallel transport along $c$\/}.

It is of some interest to rephrase what we have said at the beginning of the
discussion on $\rho$-admissible curves: if $c:I\rightarrow V$ is
$\rho$-admissible, then for every $\psi:I\rightarrow P$ which projects onto
$c_M$, the curve $t\mapsto (c(t),\dot{\psi}(t))$ in fact is a
$\rho^1$-admissible curve in $\prolr{P}$. This idea can be pushed a bit
further. Indeed, when thinking of curves in the context of our alternative
view on $\rho$-connections, it is rather the following construction which
looks like the natural thing to do.

Consider a curve $\gamma$ in $\mu^*V$, i.e.\ $\gamma$ is of the form
$\gamma:t\mapsto (\psi(t),c(t))$, with $c:I\rightarrow V$ and $\psi:I\rightarrow P$,
whereby the only assumption at the start is that $\psi_M=c_M$.
Take its horizontal lift $\H{\gamma}:I\rightarrow \prolr{P}$
which is defined, according to (\ref{Hlift}), by
\begin{equation}
t\mapsto \H{\gamma}(t)=\big(c(t),h(\psi(t),c(t))\big). \label{horcurve}
\end{equation}
Then, we could define $\psi$ to be $c^h$, the horizontal lift of $c$, if
$\H{\gamma}$ is a $\rho^1$-admissible curve in $\prolr{P}$. Indeed, it is
clear by construction that $\mu^1\circ \H{\gamma}=\psi$, so that
$\rho^1$-admissibility requires that $\dot{\psi}=\rho^1\circ
\H{\gamma}=h(\psi,c)$. Since $\psi_M=c_M$, this implies in particular that
$\dot{c}_M=T\mu\circ\dot{\psi}=T\mu(h(\psi,c))=\rho\circ c$. So, this
alternative definition implies that $c$ will necessarily have to be
$\rho$-admissible. Furthermore, from comparing what $\rho^1$-admissibility
means with (\ref{covderivcurve}) and (\ref{h-lift}), it is clear that we are
talking then about the same concept of horizontal lift $c^h$.

Note, by the way, that this other way of defining $c^h$ by no means relies on
the assumption of linearity of the $\rho$-connection. So, it is perfectly
possible to talk about parallel transport also in the context of non-linear
connections. The difference then is, of course, that if we look at points of
$P$ in the image of curves $c^h$ with different initial values in $P_m$, and
this as a map between fibres of $P$, there need not be any special feature to
talk about (compared to the fibre-wise linear action of this map we have in
the case of a linear connection); also, if $c$ has a given interval as domain,
$c^h$ need not be defined over the same domain. Needless to say, one can
introduce such a generalisation also within the more traditional approach
described first. Indeed, the map $K$ makes sense for arbitrary
$\rho$-connections and as a result one can introduce an operation
$\del{\z}\sigma$ also in this more general situation. This then still depends
on the section $\z$ of $V$ in a $\cinfty{M}$-linear way, but the fact that
such a $\del{}$ is not very commonly used comes from the failure of having a
derivation property with respect to the module structure of $Sec(\mu)$.

\section{The case of an affine bundle and its bidual}

Suppose that $\pi: E \rightarrow M$ is an affine bundle, modelled on a vector
bundle $\ov\pi: \ov E \rightarrow M$. For any $m\in M$, $E^\dagger_m :=
\mbox{Aff}(E_m,\R)$ is the set of all affine functions on $E_m$ and
$E^\dagger= \bigcup_{m\in M}E^\dagger_m$ is a vector bundle over $M$, called
the extended dual of $E$. In turn, the dual of $\pi^\dagger: E^\dagger
\rightarrow M$, denoted by $\tpi: \tilde E:= (E^\dagger)^* \rightarrow M$, is
a vector bundle into which both $E$ and $\ov E$ can be mapped via canonical
injections, denoted respectively by $\iota$ and ${\pmb \iota}$. The map
$\iota$ is affine and has ${\pmb \iota}$ as its associated linear map. With
reference to the previous section, the situation we will focus on now is the
case where $\mu:P\rightarrow M$ is the affine bundle $\pi:E\rightarrow M$,
whereas $\tau:V\rightarrow M$ still is an arbitrary vector bundle. Our main
objective is to define and characterise $\rho$-connections on $\pi$ which are
{\em affine\/}. For that purpose, we will need the overall diagram of the
previous section also with the vector bundle $\tpi:\tilde E\rightarrow M$ in
the role of $\mu:P\rightarrow M$.

\begin{dfn} A $\rho$-connection $h$ on the affine bundle $\pi:E\rightarrow M$
is said to be affine, if there exists a linear $\rho$-connection $\tilde h:
\tpi^*V\rightarrow T\tilde E$ on $\tpi:\tilde E\rightarrow M$ such that,
\[
\tilde h\circ \iota = T\iota\circ h.
\]
\end{dfn}

Both sides in the above commutative scheme of course are regarded as
maps from $\pi^*V$ to $T\tilde E$, which means that the $\iota$ on the left
stands for the obvious extension $\iota:\pi^*V\rightarrow\tpi^*V, (e,v)\mapsto
(\iota(e),v)$.

Probably the best way to see what this concept means is to look at a coordinate
representation. Let $x^i$ denote
coordinates on $M$ and $y^\a$ fibre coordinates on $E$ with respect
to some local frame $\cosyst$ for $Sec(\pi)$. The induced basis for $Sec(\pi^\dagger)$
is denoted by $(e^0,e^\alpha)$ and defined as follows: for each $a\in Sec(\pi)$ with
local representation $a(x)=e_0(x) + a^\alpha(x){\pmb e}_\alpha(x)$,
\[
e^0(a)(x)=1,\ \forall x,\qquad e^\alpha(a)(x)=a^\alpha(x).
\]
In turn, we denote the dual basis for $Sec(\tpi)$ by $(e_0,e_\alpha)$ (so that
in fact $\iota(e_0)=e_0$ and ${\pmb \iota}({\pmb e}_\alpha)=e_\alpha$).
Induced coordinates on $\tilde E$ are denoted by
$(x^i,y^A)=(x^i,y^0,y^\alpha)$. For the coordinate representation of a point
$v\in V$, we will typically write $(x^i,v^a)$. The anchor map
$\rho:V\rightarrow TM$ then takes the form $\rho:(x^i,v^a)\mapsto
\rho^i_a(x)v^a \fpd{}{x^i}$.

Following \cite{FransBavo}, we know that the map
$h:\pi^*V\rightarrow TE$ locally is of the form:
\begin{equation}
h(x^i,y^\alpha,v^a)=(x^i,y^\alpha,\rho^i_a(x)v^a,-\Gamma^\alpha_a(x,y)v^a), \label{h}
\end{equation}
whereby we have adopted a different sign convention concerning the connection
coefficients $\Gamma^\alpha_a$. Similarly, $\tilde h:\tpi^*V\rightarrow T\tilde E$,
which is further assumed to be linear, takes the form
\begin{equation}
\tilde h(x^i,y^A,v^a)=(x^i,y^A,\rho^i_a(x)v^a,-\tilde \Gamma^A_{aB}(x)y^Bv^a). \label{th}
\end{equation}
We have
\[
\tilde h(\iota(e),v)=\Big(x^i,1,y^\alpha,\rho^i_a(x)v^a,-(\tilde \Gamma^A_{a0}(x)
+\tilde \Gamma^A_{a\beta}(x)y^\beta)v^a)\Big),
\]
whereas
\[
Ti\circ h (e,v)=\Big(x^i,1,y^\alpha,\rho^i_a(x)v^a,0,-\Gamma^\alpha_a(x,y)v^a\Big).
\]
It follows that $\tilde \Gamma^0_{aB}=0$ and, more importantly, that the connection
coefficients of the affine $\rho$-connection $h$ are of the form (omitting tildes)
\begin{equation}
\Gamma^\alpha_a(x,y)=\Gamma^\alpha_{a0}(x)+\Gamma^\alpha_{a\beta}(x)y^\beta . \label{affinecon}
\end{equation}
Notice that $\ov \pi:\ov E\rightarrow M$ is a (proper) vector subbundle of
$\tpi$. With respect to the given anchor map, it of course also has its
$\rho$-prolongation $\prolr{\ov E}$. Taking the restriction of the linear
$\rho$-connection $\tilde h$ to ${\ov \pi}^*V$, we get a linear
$\rho$-connection $\ov h$ on $\ov \pi$, meaning that $\tilde{h}\circ {\pmb
\iota}=T{\pmb \iota}\circ\bar{h}$. The above coordinate expressions make this
very obvious. Indeed, if $(x^i,w^\alpha)$ are the coordinates of an element
${\pmb w}\in \ov{E}$, we have
\begin{eqnarray*}
\lefteqn{\ov{h}(x^i,w^\alpha,v^a) = \tilde{h}(x^i,0,w^\alpha,v^a) }\\
&& = (x^i,0,w^\alpha,\rho^i_av^a,0,-\Gamma^\alpha_{a\beta}w^\beta v^a) \quad
\mbox{as element of $T\tilde{E}$} \\
&& = (x^i,w^\alpha,\rho^i_av^a,-\Gamma^\alpha_{a\beta}w^\beta v^a) \quad
\mbox{as element of $T\ov{E}$}.
\end{eqnarray*}
Note further that we can formally write for the coordinate expression of
$h(e+{\pmb w},v)$:
\begin{eqnarray*}
h(x^i,y^\alpha+w^\alpha,v^a) &=& \big(x^i,y^\alpha+w^\alpha,\rho^i_av^a,
-(\Gamma^\alpha_{a0}+\Gamma^\alpha_{a\beta}y^\beta)v^a -
\Gamma^\alpha_{a\beta}w^\beta v^a\big) \\
&=& h(x^i,y^\alpha,v^a)+\ov{h}(x^i,w^\alpha,v^a).
\end{eqnarray*}
But this is more than just a formal way of writing: the following intrinsic
construction which generalises (\ref{linearcon}) is backing it. Let $\Sigma$
denote the action of $\ov E$ on $E$ which defines the affine structure, i.e.\
$\Sigma(e,{\pmb w})=e+{\pmb w}$ for $(e,{\pmb w})\in E \times_M \ov E$. Then
the above formal relation expresses that we have:
\begin{equation}
h(e+{\pmb w},v) = T_{(e,{\pmb w})}\Sigma\,\big(h(e,v),\ov h ({\pmb w}, v)\big)
\label{hhbarcond}
\end{equation}
In fact, by reading the above coordinate considerations backwards, roughly
speaking, one can see that (\ref{hhbarcond}), for a given linear $\ov h$, will
imply that the connection coefficients of the $\rho$-connection $h$ have to be
of the form (\ref{affinecon}). In other words, the following is an equivalent
definition of affineness of $h$.
\begin{dfn} A $\rho$-connection $h$ on the affine bundle $\pi:E\rightarrow M$
is affine, if there exists a linear $\rho$-connection $\ov
h:\ov{\pi}^*V\rightarrow T{\ov E}$ on $\ov \pi:\ov E\rightarrow M$, such that
(\ref{hhbarcond}) holds for all $(e,{\pmb w})\in E\times_M\ov E$.
\end{dfn}
One can then construct an extension $\tilde{h}:\tilde{\pi}^*V\rightarrow
T\tilde{E}$, which coincides with $\ov h$ when restricted to $\ov{\pi}^*V$, by
requiring that $\tilde{h}$ be linear and satisfy
$\tilde{h}\circ\iota=T\iota\circ h$.

As shown in Theorem~1, a $\rho$-connection on $\pi$ is equivalent to a decomposition of the
bundle $\prolr{E}$, originating from a horizontal lift operation from $\pi^*V$ to $\prolr{E}$ (or
sections thereof). In the representation (\ref{LVP}) of points of $\prolr{E}$ as couples of an
element of $V$ and a suitable tangent vector of $E$, the horizontal lift is given by
\[
\H{(x^i,y^\alpha,v^a)}= \left((x^i,v^a),v^a\left(\rho^i_a\fpd{}{x^i}
-\Gamma^\alpha_a\fpd{}{y^\alpha}\right)\right).
\]
At this stage, it is of interest to introduce a local basis for sections of
the $\rho$-prolongation $\pi^1:\prolr{E}\rightarrow E$. A natural choice,
adapted to the choice of a local frame in $Sec(\pi)$, the natural basis of
$\vectorfields E$ and the choice of a local basis of sections ${\pmb v}_a$ of
$\tau$, is determined as follows: for each $e\in E$, if $x$ are the
coordinates of $\pi(e)\in M$,
\begin{equation}
\baseX a (e)=
\left({\pmb v}_a(x),\left.\r^i_{a}(x)\fpd{}{x^i} \right|_e\right), \quad
\baseV \a(e) = \left(0,\left.\fpd{}{y^{\a}}\right|_e\right). \label{prolongedbasis}
\end{equation}
Coordinates of a point $(v,X_e)\in\prolr{E}$ are of the form:
$(x^i,y^\alpha,v^a,X^\alpha)$. A general section of the
$\rho$-prolongation can be represented locally in the form:
\begin{equation}
{\mathcal Z}= \z^a(x,y)\baseX{a} + Z^\a(x,y)\baseV{\a}.
\label{prolZ}
\end{equation}
Its projection onto $Sec(p)$ ($p:\pi^*V\rightarrow E$) is
$\z=\z^a{\pmb v}_a$. Now, once we have a given $\rho$-connection
on $\pi$ (affine or not), we are led to introduce a local basis
for the horizontal sections of $\pi^1$, which is given by
\begin{equation}
\baseH{a}=\subH{P}(\baseX{a})=\baseX{a} - \G_a^\alpha(x,y) \baseV \a.
\label{hor}
\end{equation}
A better representation of the section (\ref{prolZ}), adapted to the given connection, then
becomes:
\begin{equation}
{\mathcal Z}= \z^\a\baseH{a} +
(Z^\a+\Gamma^\alpha_b \z^b) \baseV{\a}.
\label{prolZ2}
\end{equation}
Let us repeat that, as a result of Theorem~1 and Definition~4, the
existence of an affine $\rho$-connection on $\pi$ is equivalent to
the existence of a horizontal lift from $Sec(p)$ to $Sec(\pi^1)$,
giving rise to a direct sum decomposition (\ref{directsum}), and
which is such that, in coordinates, the connection coefficients
(\ref{hor}) are of the form (\ref{affinecon}).

We next turn our attention to the concept of connection map, and want to see for the particular
case of an affine $\rho$-connection, to what extent it gives rise also to a covariant
derivative operator and a notion of parallel transport.

When considering the $\rho$-prolongation of different bundles $P$, it is
convenient to indicate the dependence on $P$ also in the map $\rho^1$. Given a
$\rho$-connection $h$ on the affine bundle $\pi:E\rightarrow M$, the map
$\rho^1_E-h\circ j:\prolr{E}\rightarrow TE$ gives rise (as before) to a
vertical tangent vector to $E$, at the point $e$ say. As such, this vector can
be identified with an element of $\ov E$, the vector bundle on which $E$ is
modelled, at the point $\pi(e)$. With the same notational simplification as
before, we thus get a connection map
\begin{equation}
K:=\rho^1_E-h\circ j:\prolr{E}\rightarrow \ov{E}. \label{conmap1}
\end{equation}
$K$ of course also extends to a map from $Sec(\pi^1)$ to $Sec(\ov
\pi)$. It follows directly from the definition that we have
\begin{equation}
K(\baseH{a})=0, \qquad K(\baseV{\alpha})=\base{\alpha}. \label{Konsections}
\end{equation}
We wish to come back here in some more detail to the relation between the map
$K$ and the vertical projector $\subV{P}=id-\subH{P}$, coming from the direct
sum decomposition of $\prolr{E}$. In the present case of an affine bundle
$\pi:E\rightarrow M$ over a vector bundle $\ov{\pi}:\ov{E}\rightarrow M$,
there is a natural vertical lift operation from $\ov{E}_m$ to $T_eE$ for each
$e\in E_m$. It is determined by: ${\pmb w}\mapsto \V{w}_e$, where for each
$f\in\cinfty{E}$,
\[
\V{w}_e(f) = \left. \frac{d}{dt} f(e+t{\pmb w}) \right|_{t=0}.
\]
This in turn extends to an operator $\V{}:\pi^*\ov E\rightarrow \prolr{E}$,
determined by $\V{(e,{\pmb w})}=(0,\V{w}_e)$, which defines an isomorphism
between $\pi^*\ov E$ and ${\rm Im}\,\subV{P}$. The short exact sequence
(\ref{shortexactgeneral}) of which a $\rho$-connection is a splitting, can
thus be replaced by
\begin{equation}
0 \rightarrow \pi^*\ov{E} \stackrel{V}{\rightarrow} \prolr{E}
\stackrel{j}{\rightarrow} \mu^* V \rightarrow 0. \label{newshortexact}
\end{equation}
Within this picture of $\rho$-connections, the connection map $K$ thus is
essentially the co-splitting of the splitting $\H{}$, that is to say, we have
$K\circ\V{}=id_{\pi^*\ov{E}}\ $ and $\V{}\circ K + \H{}\circ
j=id_{\prolr{E}}$.

The map $K$ becomes more interesting when the connection is affine. Indeed, denoting the
projection of $\prolr{\tilde E}$ onto $\tpi^*V$ by $\tilde{j}$, it then follows from
Definition~4 that we also have a connection map
\begin{equation}
\tilde K:=\rho^1_{\tilde E}-\tilde h\circ\tilde j:\prolr{\tilde E}\rightarrow
\tilde{E}. \label{conmap2}
\end{equation}
The map $T\iota:TE\rightarrow T\tilde E$ extends to a map from $\prolr{E}$ to
$\prolr{\tilde E}$ in the following obvious way: $T\iota:(v,X_e)\mapsto (v,T\iota(X_e))$.
Indeed, we have $T\tpi(T\iota(X_e))=T(\tpi\circ i)(X_e)=T\pi(X_e)=\rho(v)$, as required.

\begin{prop}
For an affine $\rho$-connection on $\pi$ we have
\begin{equation}
{\pmb \iota}\circ K=\tilde K\circ T\iota. \label{KtildeK}
\end{equation}
\end{prop}

\proof\ \ In coordinates, $K$ and $\tilde K$ are given by
\begin{eqnarray*}
K: && (x^i,v^a,y^\alpha,Z^\alpha)\mapsto (Z^\alpha +\Gamma^\alpha_a v^a)\,\base{\alpha}(x) \\
\tilde K: && (x^i,v^a,y^A,Z^A)\mapsto Z^0 e_0(x) + (Z^\alpha +\Gamma^\alpha_{aB}y^B v^a)\,e_\alpha(x).
\end{eqnarray*}
Hence,
\begin{eqnarray*}
\tilde K\circ T\iota(x^i,v^a,y^\alpha,Z^\alpha)&=& \tilde K(x^i,v^a,1,y^\alpha,0,Z^\alpha) \\
&=& \Big(Z^\alpha +(\Gamma^\alpha_{a0}+\Gamma^\alpha_{a\beta}y^\beta) v^a\Big)\,e_\alpha(x),
\end{eqnarray*}
from which the result follows in view of (\ref{affinecon}). \qed

Notice that $\ov h$ also has a corresponding connection map $\ov K: \prolr{\ov
E}\rightarrow \ov E$, which obviously coincides with $\tilde K|_{\prolr{\ov
E}}$, so that we also have
\begin{equation}
{\pmb \iota}\circ \ov K=\tilde K\circ T{\pmb \iota}. \label{ovKtildeK}
\end{equation}

Let now $\z$ be a section of $\tau$ and $\sigma$ a section of $\pi$. If we apply
the tangent map $T\sigma:TM\rightarrow TE$ to $\rho(\z(m))$, it is obvious by
construction that $\big(\z(m),T\sigma(\rho(\z(m)))\big)$ will be an element of $\prolr{E}$.
The connection map $K$ maps this into a point of $\ov E|_{m}$. Hence, the covariant derivative
operator of interest in this context is the map $\del{}:Sec(\tau)\times Sec(\pi)
\rightarrow Sec(\ov \pi)$, defined by
\begin{equation}
\del{\z}\sigma (m)=K\big(\z(m),T\sigma(\rho(\z(m)))\big). \label{del}
\end{equation}
To discover the properties which uniquely characterise the covariant derivative associated
to an affine $\rho$-connection, we merely have to exploit the results of Proposition~2.
In doing so, we will of course rely on the known properties (see \cite{FransBavo}) of the
covariant derivative $\tilde{\del{}}$, associated to the linear $\rho$-connection $\tilde h$.
We observe that $\del{}$ is manifestly $\R$-linear in its first argument and now further
look at its behaviour with respect to the $\cinfty{M}$-module structure
on $Sec(\tau)$. From (\ref{KtildeK}), it follows that for $f\in\cinfty{M}$,
\begin{eqnarray*}
{\pmb \iota}\big((\del{f\z}\sigma)(m)\big) &=& {\pmb \iota}\Big(K\big(f\z(m),T\sigma(\rho(f\z(m)))\big)\Big) \\
&=& \tilde K\big(f\z(m),T(\iota\sigma)(\rho(f\z(m)))\big) \\
&=& \tilde{\nabla}_{f\z}(\iota\sigma)(m) = f(m)\,\tilde{\nabla}_\z(\iota\sigma)(m) \\
&=& f(m) \tilde K\big(\z(m),T\iota\circ T\sigma(\rho(\z(m)))\big) \\
&=& f(m)\, {\pmb \iota}\Big(K\big(\z(m),T\sigma(\rho(\z(m)))\big)\Big) \\
&=& {\pmb \iota}\big(f(m)\,\del{\z}\sigma(m)\big),
\end{eqnarray*}
from which it follows that
\begin{equation}
\del{f\z}\sigma = f\,\del{\z}\sigma. \label{delf1}
\end{equation}
For the behaviour in the second argument, we replace $\sigma$ by $\sigma+f{\pmb \eta}$, with
$f\in\cinfty{M}$ and ${\pmb \eta}\in Sec(\ov \pi)$.
Denoting the linear covariant derivative coming from the restriction $\ov K$ by $\ov{\nabla}$,
we compute in the same way, using (\ref{KtildeK}) and (\ref{ovKtildeK}):
\begin{eqnarray*}
\lefteqn{
{\pmb \iota}\big(\del{\z}(\sigma+f{\pmb \eta})(m)\big) = {\pmb \iota} \Big(K\big(\z(m),
T(\sigma+f{\pmb \eta})(\rho(\z(m)))\big)\Big) }\\
&& = \tilde K\big(\z(m),T(\iota\sigma+f\pmb{\iota}\pmb{\eta})(\rho(\z(m)))\big)
= \tilde{\nabla}_{\z}(\iota\sigma + f\pmb{\iota}\pmb{\eta})(m) \\
&&= \tilde{\nabla}_{\z}\iota\sigma(m) +f(m)\big(\tilde{\nabla}_\z\pmb{\iota} \pmb{\eta}\big)(m)
+ \rho(\z)(f)(m)\,\pmb{\iota}\pmb{\eta}(m)\\
&&= \tilde K\big(\z(m),T\iota\circ T\sigma(\rho(\z(m)))\big)
+ f(m)\,\tilde K\big(\z(m),T{\pmb \iota}\circ T{\pmb \eta}(\rho(\z(m)))\big) \\
&& \mbox{\ \ } + \rho(\z)(f)(m)\,\pmb{\iota}\pmb{\eta}(m)
\;= \;{\pmb \iota}\Big(K\big(\z(m),T\sigma(\rho(\z(m)))\big)\Big) \\
&& \mbox{\ \ }+ f(m)\,{\pmb \iota}\Big(\ov K\big(\z(m),T{\pmb \eta}(\rho(\z(m))\big)\Big)
+ \rho(\z)(f)(m)\,\pmb{\iota}\pmb{\eta}(m)\\
&&= {\pmb \iota}\Big(\del{\z}\sigma(m)+ f(m)\,\ov{\nabla}_\z{\pmb \eta}(m) +
\rho(\z)(f)(m){\pmb \eta}(m)\Big).
\end{eqnarray*}
This expresses that we have the property:
\begin{equation}
\del{\z}(\sigma+f{\pmb \eta})= \del{\z}\sigma + f\,\ov{\nabla}_\z{\pmb \eta} + \rho(\z)(f)\,
{\pmb \eta}. \label{delf2}
\end{equation}

In coordinates we have, for $\z=\zeta^a(x){\pmb v}_a$ and
$\sigma=e_0+\sigma^\alpha(x)\base{\alpha}$:
\begin{equation}
\del{\z}\sigma= \left(\fpd{\sigma^\alpha}{x^i}\rho^i_a(x) +
\Gamma^\alpha_{a0}(x)+\Gamma^\alpha_{a\beta}(x)\sigma^\beta(x)\right)\z^a(x)
\, \base{\alpha}. \label{27}
\end{equation}
As one can see, the linearity in $\z$ makes that the value of $\del{\z}\sigma$
at a point $m$ only depends of the value of $\z$ at $m$, so that the usual
extension works, whereby for any fixed $v\in V$, $\del{v}$ is a map from
$Sec(\mu)$ to $\ov{E}_m$, defined by $\del{v}\sigma = \del{\z}\sigma (m)$, for
any $\z$ such that $\z(m)=v$.

\begin{thm}
An affine $\rho$-connection $h$ on $\pi$ is uniquely characterised by the existence of an
operator $\del{}:Sec(\tau)\times Sec(\pi) \rightarrow Sec(\ov \pi)$ and an associated
$\ov{\nabla}:Sec(\tau)\times Sec(\ov \pi) \rightarrow Sec(\ov \pi)$, such that $\del{}$
is $\R$-linear in its first argument, $\ov{\nabla}$ satisfies the requirements for the
determination of a linear $\rho$-connection on $\ov \pi$, and the properties (\ref{delf1})
and (\ref{delf2}) hold true.
\end{thm}

\proof\ \ Given an affine $\rho$-connection $h$ on $\pi$, the existence of
operators $\del{}$ and $\ov{\nabla}$ with the required properties has been
demonstrated above. Assume conversely that such operators are given. Then,
there exists an extension $\tilde{\nabla}:Sec(\tau)\times Sec(\tpi)
\rightarrow Sec(\tpi)$, which is defined as follows. Every $\tilde{\sigma} \in
Sec(\tpi)$ locally is either of the form $\tilde{\sigma}=f\,\iota(\sigma)$ for
some $\sigma\in Sec(\pi)$ or of the form $\tilde{\sigma}={\pmb \iota}({\pmb
\eta})$ for some ${\pmb \eta}\in Sec(\ov \pi)$. In the first case, we put
\[
\tilde{\nabla}_\z\tilde{\sigma} = f\,{\pmb \iota}(\del{\z}\sigma) +
\rho(\z)(f)\iota(\sigma);
\]
in the second case, we put
\[
\tilde{\nabla}_\z\tilde{\sigma} = {\pmb \iota}(\ov{\nabla}_\z{\pmb \eta}).
\]
We further impose $\tilde{\nabla}$ to be $\R$-linear in its second argument.
$\R$-linearity as well as $\cinfty{M}$-linearity in the first argument
trivially follows from the construction. It is further easy to verify that for
$g\in\cinfty{M}$: $\tilde{\nabla}_\z(g \tilde{\sigma}) =
g\,\tilde{\nabla}_\z\tilde{\sigma} + \rho(\z)(g)\,\tilde{\sigma}$. Indeed, in
the case that $\tilde{\sigma}=f\,\iota(\sigma)$, for example, we have
\begin{eqnarray*}
\tilde{\nabla}_\z(g \tilde{\sigma}) &=& gf\,{\pmb \iota}(\del{\z}\sigma) +
\big( f\,\rho(\z)(g) + g\,\rho(\z)(f)\big)\iota(\sigma) \\
&=& g\,\tilde{\nabla}_\z\tilde{\sigma} + \rho(\z)(g)\,\tilde{\sigma},
\end{eqnarray*}
and likewise for the other case. Following \cite{FransBavo} we thus conclude
that $\tilde{\nabla}$ uniquely determines a linear $\rho$-connection on $\tpi$
by the following construction: for each $(\tilde{e},v)\in\tpi^*V$, take any
$\tilde \psi\in Sec(\tpi)$ for which $\tilde \psi(\tau(v))=\tilde{e}$, and put
\[
\tilde{h}(\tilde{e},v)= T\tilde{\psi}(\rho(v)) -
\V{(\tilde{\nabla}_v\tilde{\psi})}_{\tilde e},
\]
where the last term stands for the element
$\tilde{\nabla}_v\tilde{\psi}(\tau(v))\in \tilde{E}_{\tau(v)}$, vertically
lifted to a vector tangent to the fibre of $\tilde{E}$ at $\tilde{e}$.

Likewise, we define a fibre linear map $h:\pi^*V\rightarrow TE$ by
\[
h(e,v)= T\psi(\rho(v)) - \V{(\nabla_v\psi)}_{e},
\]
which can be seen to be independent of the choice of a section $\psi$ for
which $\psi(\tau(v))=e$. It is obvious that $h$ satisfies the requirements of
a $\rho$-connection on $\pi$. It remains to show that
$\tilde{h}\circ\iota=T\iota\circ h$. We have
\begin{eqnarray*}
\tilde{h}(\iota(e),v)&=& T(\iota\psi)(\rho(v)) -
\V{\big(\tilde{\nabla}_v(\iota\psi)\big)}_{\iota(e)} \\
&=& T(\iota\psi)(\rho(v)) - \V{({\pmb \iota}\del{v}\psi)}_{\iota(e)} \\
&=& T\iota\circ T\psi(\rho(v)) - T\iota\big(\V{(\del{v}\psi)}_e\big) \\
&=& T\iota(h(e,v)),
\end{eqnarray*}
which completes the proof. \qed

Another interesting question one can raise in this context is about the
circumstances under which a linear $\rho$-connection $\tilde{h}$ on
$\tilde{\pi}$ is associated to an affine $\rho$-connection $h$ on $\pi$ in the
sense of Definition~4. A simple look at coordinate expressions leads to the
following result with a global meaning.

\begin{prop} A linear $\rho$-connection on
$\tilde{\pi}$ is associated to an affine $\rho$-connection on $\pi$ if and
only if $e^0$ is parallel.
\end{prop}

\proof\ \ For the covariant derivative operator $\tilde{\nabla}$ associated to
a linear $\tilde{h}$, we have for the local basis of $Sec(\tilde{\pi})$:
\[
\tilde{\nabla}_\zeta e_A = \zeta^a \tilde{\Gamma}^B_{aA}\,e_B,
\]
and by duality, for the basis of $Sec(\pi^\dagger)$:
\[
\tilde{\nabla}_\zeta e^A = -\zeta^a \tilde{\Gamma}^A_{aB}\,e^B.
\]
It follows that $\tilde{\nabla}_\zeta e^0=0 \ \Leftrightarrow \
\tilde{\Gamma}^0_{aB}=0$. The restriction of $\tilde{h}$ to $\iota(E)$ then
defines an affine $\rho$-connection on $\pi$. \qed

A few words are in order, finally,  about the concept of parallel transport in
this case. Following the comments about $\rho^1$-admissibility of a curve
$\H{\gamma}$ made at the end of the previous section, we know that a curve
$\psi$ in $E$, with coordinate representation $t\mapsto
(x^i(t),\psi^\alpha(t))$ will be the horizontal lift $c^h$ of a
$\rho$-admissible curve $c:t\mapsto (x^i(t),c^a(t))$, provided that (cf. the
coordinate expressions (\ref{h}) and (\ref{affinecon})) $x^i(t)$ and
$\psi^\alpha(t)$ satisfy the differential equations:
\begin{eqnarray}
\dot{x}^i &=& \rho^i_a(x)c^a(t), \label{radm} \\
\dot{\psi}^\alpha &=& - \Gamma^\alpha_{a0}(x)c^a(t)
- \Gamma^\alpha_{a\beta}(x)c^a(t)\,\psi^\beta. \label{partp}
\end{eqnarray}
In the more standard approach to the definition of $c^h$, if $c$ is a
$\rho$-admissible curve in $V$ and $\psi$ a curve in $E$ which projects onto
$c_M$, we can define a new curve $\del{c}\psi$ by a formula which is formally
identical to (\ref{covderivcurve}). Note, however, that $\del{c}\psi$ is a
curve in $\ov E$ now. Nevertheless, it makes perfect sense to say that $\psi$
in $E$ is $c^h$ if the associated curve $\del{c}\psi$ in $\ov E$ is zero for
all $t$. It can be seen from the coordinate expression (\ref{partp}) that for
different initial values in a fixed fibre of $E$, we get an affine action
between the affine fibres of $E$, whose corresponding linear part comes from
the parallel transport rule associated to the linear connection $\ov h$ on
$\ov E$. This is in agreement with the property, coming from (\ref{delf2}),
that
\begin{equation}
\del{\z}(\sigma+{\pmb \eta})= \del{\z}\sigma + \ov{\nabla}_\z{\pmb \eta}.
\end{equation}

\section{Discussion}

Perhaps the simplest example of the natural appearance of an affine $\rho$-connection
(though for a trivial $\rho$), is the following. Take $E$ to be the first-jet bundle
$J^1M$ of a manifold $M$ which is fibred over $\R$, and $V=TM$ with $\rho=id_{TM}$.
Then $\prolr{E}=TE$ and we are in the situation which has been extensively studied
in \cite{CraMaSa}. It is well-known that every second-order differential equation
field (\sode) on $J^1M$, say
\[
\Gamma= \fpd{}{t}+ v^i\fpd{}{x^i} + f^i(t,x,v)\fpd{}{v^i},
\]
defines a non-linear connection whose connection coefficients are
\[
\Gamma^i_j= - \frac{1}{2}\fpd{f^i}{v^j} \qquad \G^i_0= - f^i
+\frac{1}{2}\fpd{f^i}{v^j}v^j.
\]
To say that the forces $f^i$ are quadratic in the velocities, i.e. are of the form
\[
f^i = f^i_0(t,x) + f^i_j(t,x)v^j + f^i_{jk}(t,x)v^jv^k,
\]
is an invariant condition and clearly gives rise then to a connection of affine type,
as discussed in the previous section.

For this standard example, however, there is more structure available then
merely this connection and that is what makes the geometrical study of \sode s such
a rich subject. The extra structure primarily comes from two sides. First of all,
there is the structure of $J^1M$ itself where, in particular, a canonical vertical
endomorphism is defined (which in fact lies at the origin of the \sode-connection,
see e.g. \cite{Crampin}). Secondly, since sections of $V$ and of $\prolr{E}$ here
are simply vector fields, they come equipped with a Lie algebra structure and this
in turn is essential for defining such concepts as torsion and curvature of a
connection.

We intend to study in a forthcoming paper a quite general situation of affine
$\rho$-connections, where the same kind of extra structure is available. To that end
we will take $\pi:E\rightarrow M$ to be a general affine bundle and let the
vector bundle $\tau:V\rightarrow M$ be the bidual $\tpi:\tilde E\rightarrow M$.
In addition we will assume that $E$ comes equipped with an affine Lie algebroid
structure (as studied for example in \cite{affalg1,affalg2}). This implies that
the anchor map $\rho:V\rightarrow TM$ then also becomes the anchor of a (vector)
Lie algebroid. As shown in \cite{affalg2}, the prolonged bundle $\prolr{E}$
inherits a Lie algebroid structure; moreover there is a canonical endomorphism
on sections of $\prolr{E}$, which is exactly the analogue of the vertical
endomorphism on a first-jet bundle. Not surprisingly therefore, it is possible
to define dynamical systems of Lagrangian type on such an affine Lie algebroid.
Much of this has been explored already in the above cited papers, but the theory of
affine connections and so-called pseudo-\sode s in that context still needs
to be developed.

{\small {\bf Acknowledgements.} E.\ Mart\'{\i}nez acknowledges partial
financial support from CICYT grant BFM2000-1066-C03-01.

}

\end{document}